\documentclass[11pt, draft]{amsart}
\usepackage{amssymb}

\theoremstyle{plain}
\newtheorem{theorem}{Theorem}
\newtheorem{lemma}{Lemma}

\theoremstyle{remark}
\newtheorem{remark}{Remark}

\begin{document}
\title[Preserving maximal deviation of observables]{Linear maps
on the space of all bounded observables preserving maximal deviation}
\author{LAJOS MOLN\'AR}
\address{Institute of Mathematics and Informatics\\
         University of Debrecen\\
         4010 Debrecen, P.O.Box 12, Hungary}
\email{molnarl@math.klte.hu}
\author{M\'ATY\'AS BARCZY}
\address{Institute of Mathematics and Informatics\\
         University of Debrecen\\
         4010 Debrecen, P.O.Box 12, Hungary}
\email{barczy@math.klte.hu}
\thanks{This paper was written when the first author held a Humboldt
Research Fellowship. He is very grateful to
the Alexander von Humboldt Foundation for providing ideal conditions
for research and also to his host Werner Timmermann (TU Dresden, Germany)
for very warm hospitality.
         The first author also acknowledges support from
         the Hungarian National Foundation for Scientific Research
         (OTKA), Grant No. T030082, T031995, and from
         the Ministry of Education, Hungary, Grant
         No. FKFP 0349/2000}
\date{{May 7, 2003}}
\begin{abstract}
In this paper we determine all the bijective linear maps on the
space of bounded observables which preserve a fixed moment or
the variance. Nonlinear versions of the corresponding results
are also presented.
\end{abstract}
\maketitle

\section{Introduction and Statements of the Results}

In the Hilbert space formalism of quantum mechanics there are several
structures of linear operators which play distinguished role in the
theory. These are, among others, the following. The Jordan
algebra $B_s(H)$ of all self-adjoint bounded linear operators on the
Hilbert space $H$
which are called bounded observables, the lattice $P(H)$
of all projections (i.e., self-adjoint idempotents) on $H$ called
quantum events, the convex
set $S(H)$ of all positive trace-class operators on $H$ with trace 1
called (mixed) states, and the so-called effect algebra
$E(H)$ of all positive bounded linear operators which are majorized by
the identity $I$. These structures play essential role in the
probabilistic aspects of quantum theory.

Just as in the case of any algebraic structure in mathematics in
general,
the study of the automorphisms of the above mentioned structures is of
remarkable importance. One can find an interesting unified treatment of
those automorphisms in \cite{CVLL}.
In our recent papers \cite{Mo1, BaTo} we presented some
results on the local behaviour of the automorphisms in question, while
in \cite{Mo2, Mo3, Mo4} we have started to study how these automorphisms
can be characterized by their preserving properties.

The systematic study of preserver problems (more precisely, linear
preserver problems, so-called LPP's) constitutes a part of matrix theory.
In fact, this area represents one of the most active research fields in
matrix theory (we refer only to two survey papers \cite{LiTs, LiPi}).
In the last decades considerable attention has also been paid to
the infinite dimensional case as well, i.e., to linear preserver
problems concerning algebras of linear operators on general Hilbert
spaces or Banach spaces (once again, we only refer to a survey
paper \cite{BrSe3}).
From the point of view of the present paper,
the most important point is that
the solutions of linear preserver problems provide, in most of the
cases, important new information on the automorphisms
of the underlying algebras (matrix algebras, or more generally,
operator algebras)
as they show how those automorphisms
are determined by their various preserving properties. These
properties mainly concern a
certain important numerical quantity or a set of them corresponding to
operators (e.g., norm, spectrum), or they concern a distinguished set of
operators (e.g., the set of projections),
or they concern an important relation among operators (e.g.,
commutativity).
This kind of results may help to better understand the behaviour of the
automorphisms of the underlying algebras.

In our above mentioned papers \cite{Mo2, Mo3, Mo4} we have started to
study the automorphisms of Hilbert space effect algebras and those of
the Jordan algebra of bounded observables from a similar, preserver
point of view. There we have considered transformations which preserve
quantities, or relations, or properties that all have physical meaning.
For example, as for observables, in \cite{Mo3} we determined all
bijective transformations (no linearity was assumed) of $B_s(H)$
that preserve the order (which is just the usual order among
self-adjoint operators).
In \cite{Mo4} we described the general form of those bijections of
$B_s(H)$ which preserve commutativity (in quantum theory the expression
compatibility is used in the place of commutativity) and are
multiplicative on commuting pairs of operators.

We now turn to the content of the present paper. In classical
probability theory the mean value (or, more generally, the moments) and
the variance are among the most important characteristics of a random
variable. Therefore, it
is not surprising that the same is true for the quantum mechanical
variables, i.e., for the observables. The main aim of this paper
is to show that the preservation of any of those quantities
more or less completely characterizes the automorphisms among the linear
transformations of $B_s(H)$.

In what follows, let $H$ be a complex Hilbert space.
Let $A\in B_s(H)$ and pick a unit vector $\varphi\in H$.
The mean value $m(A,\varphi)$ of the
observable $A$ in the (pure) state represented by $\varphi$ is
defined as
\[
m(A,\varphi)=\langle A\varphi, \varphi\rangle.
\]
So, unlike in classical probability, in quantum theory there is a set of
mean values of a single variable. We intend to determine all the
bijective linear transformations $\phi$ of $B_s(H)$ which preserve this
set in the sense that
\begin{align*}
\{ m(\phi(A),\varphi) \, : \, \varphi \in H, \| \varphi\|=1\} &=
\{ \langle \phi(A)\varphi,\varphi\rangle \, : \, \varphi \in H, \|
\varphi\|=1\} \\
& = \{ \langle A\varphi,\varphi\rangle \, : \, \varphi \in H, \|
\varphi\|=1\} \\
& = \{ m(A,\varphi) \, : \, \varphi \in H, \| \varphi\|=1\} \\
\end{align*}
holds for every $A\in B_s(H)$.
Clearly, the set of all mean values of an observable $A\in B_s(H)$ is
equal to the numerical range of the operator $A$.
So, the above problem can be reformulated as the linear preserver
problem concerning the numerical range on $B_s(H)$. Obviously,
it is a more general problem to preserve the numerical radius $w(.)$
instead of the numerical range. It is well-known that for
a self-adjoint operator $A$ this former quantity $w(A)$
is equal to the operator norm $\| A\|$.
Hence, we easily arrive at the problem of describing the
surjective linear isometries of $B_s(H)$. The solution of this
problem is well-known in the literature.
For example, one can consult the paper \cite{DFR}.
The corresponding result reads as follows.

\begin{theorem}\label{T:variance1}
Let $\phi:B_s(H) \to B_s(H)$ be a bijective linear map which preserves
the operator norm, that is, suppose that
\[
\| \phi(A)\|=\|A\| \qquad (A\in B_s(H)).
\]
Then there is an either unitary or
antiunitary operator $U$ on $H$ such that $\phi$ is either of the form
\begin{equation}\label{E:var11}
\phi(A)=UAU^* \qquad (A\in B_s(H))
\end{equation}
or of the form
\begin{equation}\label{E:var12}
\phi(A)=-UAU^* \qquad (A\in B_s(H)).
\end{equation}
\end{theorem}

(By an antiunitary operator we mean a norm preserving conjugate-linear
bijection of the underlying Hilbert space $H$.)
Although this is not a new result, in Section 2 we present the sketch of
a short proof that applies preserver techniques.

Observe that the above statement is a self-adjoint analogue of a
well-know result of Kadison \cite{Kad} on the surjective isometries of
$C^*$-algebras and also that of a result of Bre\v sar and \v Semrl
\cite{BrSe3} describing the form of all bijective linear maps of the
algebra of all bounded linear operators on a Banach space
which preserve
the spectral radius (recall that the norm of a self-adjoint operator is
equal to its spectral radius). However, there is no doubt,
those results are much deeper than the one we have formulated above.

By the help of Theorem~\ref{T:variance1} we can describe the bijective
linear
maps of $B_s(H)$ which preserve the set of mean values. In fact, as the
second possibility \eqref{E:var12} can be excluded, we obtain that the
maps in question are exactly the automorphisms of the Jordan algebra
$B_s(H)$ (cf. \cite{CVLL}). Moreover, observe that using the same result
Theorem~\ref{T:variance1} we can solve also the problem of preserving a
fixed moment of bounded observables. For any $n\in \mathbb N$, the
$n$th moment of an observable $A\in B_s(H)$ is
the set
\[ \{ m(A^n,\varphi) \, :\, \varphi \in H, \|\varphi\|=1\}=
\{ \langle A^n\varphi,\varphi \rangle \, :\, \varphi \in H,
\|\varphi\|=1\}.
\]
Now, the solution of the mentioned problem immediately follows as one
can refer to the equality
\[
\sup \{ | \langle A^n\varphi,\varphi \rangle | \, :\, \varphi \in H,
\|\varphi\|=1\}=
w(A^n)=\|A^n\|=\|A\|^n
\]
which holds for every self-adjoint operator $A$ on $H$.

Beside moments, the other very important probabilistic character of an
observable is its variance. Just as with mean values, we have
variance with respect to every (pure) state.
Let $A\in B_s(H)$ and $\varphi \in H,
\|\varphi\|=1$. The variance $var(A,\varphi)$ of $A$ in the state
$\varphi$ is defined by
\begin{align*}
var(A,\varphi)&= m((A-m(A,\varphi)I)^2,\varphi)\\
&= \langle(A-\langle A\varphi,\varphi\rangle
   I)^2\varphi,\varphi\rangle\\
&= \langle A^2\varphi,\varphi\rangle-\langle A\varphi,\varphi\rangle^2.
\end{align*}

We intend to determine all bijective linear maps on $B_s(H)$
which preserve the set of variances of observables.
It is obvious
that every linear map $\phi$ on $B_s(H)$ which preserves this set, i.e.,
which satisfies
\[
\{ var(\phi(A),\varphi) \, :\, \varphi\in H, \| \varphi\|=1\}=
\{ var(A,\varphi) \, :\, \varphi\in H, \| \varphi\|=1\}
\]
for every $A\in B_s(H)$, also preserves the quantity
\begin{equation}\label{E:var51}
\|A\|_v=\sup_{ \|\varphi\|=1} var(A,\varphi)^{1/2},
\end{equation}
i.e.,
satisfies
\[
\| \phi(A)\|_v=\| A\|_v
\]
for every $A\in B_s(H)$.
The quantity $\| A\|_v$ is called the maximal deviation of the
observable $A\in B_s(H)$. In its definition \eqref{E:var51} we have used
the square root
of the variances since, as it will be clear from Lemma~\ref{L:var1},
the so-obtained quantity is a semi-norm
on $B_s(H)$ which is quite convenient to handle.

Observe that every automorphism of $B_s(H)$ (see \cite{CVLL}) as well as
its negative preserves the
maximal deviation and that perturbations by scalar operators also do not
change this quantity. Our result that follows (which can be considered
as the main result of the paper) states that from these two types of
transformations we can construct all the linear preservers under
consideration.

\begin{theorem}\label{T:variance2}
Let $\phi:B_s(H)\to B_s(H)$ be a bijective linear map which preserves
the maximal deviation, that is, suppose that
\[
\|\phi(A)\|_v=\|A\|_v \qquad (A\in B_s(H)).
\]
Then there exist an either unitary or
antiunitary operator $U$ on $H$ and a linear functional
$f:B_s(H)\to\mathbb R$ such that $\phi$ is either of the form
\begin{equation}\label{E:var13}
\phi(A)=U AU^*+f(A)I \qquad (A\in B_s(H))
\end{equation}
or of the form
\begin{equation}\label{E:var14}
\phi(A)=-U AU^*+f(A)I \qquad (A\in B_s(H)).
\end{equation}
\end{theorem}

Unlike with the transformations preserving the set of mean values,
for the bijective linear maps on $B_s(H)$ which preserve the set of
variances, the second possibility \eqref{E:var14} above can obviously
occur. Hence, we obtain that every such preserver is "an automorphism of
$B_s(H)$ or its negative perturbed by a scalar operator valued linear
transformation".

Since, from the physical point of view, to assume the linearity
of the considered transformations on the space of observables
sometimes seems to be a strong assumption that can be quite difficult to
check in the particular cases, in the remaining results
we formulate nonlinear versions of Theorems~\ref{T:variance1} and
\ref{T:variance2} as follows.
First observe that
\[
d_m(A,B)=\sup_{ \|\varphi\|=1} |m(A-B,\varphi)|=\| A-B\|
\quad (A,B\in B_s(H))
\]
defines a metric on $B_s(H)$, while
\[
d_v(A,B)=\sup_{ \|\varphi\|=1} var(A-B,\varphi)^{1/2}=\| A-B\|_v
\quad (A,B\in B_s(H))
\]
defines a semi-metric on $B_s(H)$. Both $d_m$ and $d_v$
represent certain stochastic distances between bounded observables.
Using the first two results and the celebrated Mazur-Ulam theorem
on surjective nonlinear isometries of normed spaces \cite{MaUl}, we can
prove the following statements
which show how close the stochastic isometries with respect to
either $d_m$
or $d_v$ are to the automorphisms of the Jordan algebra
$B_s(H)$.

\begin{theorem}\label{T:variance3}
Let $\phi:B_s(H) \to B_s(H)$ be a bijective transformation (linearity is
not assumed) with the property that
\[
d_m(\phi(A),\phi(B))=d_m(A,B) \qquad (A\in B_s(H)).
\]
Then there are an either unitary or
antiunitary operator $U$ on $H$ and a fixed operator $X\in B_s(H)$ such
that $\phi$ is either of the form
\[
\phi(A)=UAU^* +X \qquad (A\in B_s(H))
\]
or of the form
\[
\phi(A)=-UAU^* +X \qquad (A\in B_s(H)).
\]
\end{theorem}

The last result of the paper describes the form of all
"stochastic isometries" with respect to the semi-metric $d_v$.

\begin{theorem}\label{T:variance4}
Let $\phi:B_s(H)\to B_s(H)$ be a bijective transformation (linearity is
not assumed) with the property that
\[
d_v(\phi(A),\phi(B))=d_v(A,B) \qquad (A\in B_s(H)).
\]
Then there exist an either unitary or
antiunitary operator $U$ on $H$, a fixed operator $X\in B_s(H)$, and a
functional
$f:B_s(H)\to\mathbb R$ (not linear in general) such that $\phi$ is
either of the form
\[
\phi(A)=U AU^*+X+f(A)I \qquad (A\in B_s(H))
\]
or of the form
\[
\phi(A)=-U AU^*+X+f(A)I \qquad (A\in B_s(H)).
\]
\end{theorem}

\section{Proofs}

We first remark that in what follows whenever we speak about the
preservation of an object or relation we always mean that this is
preserved in both directions.

We now present a short proof of Theorem~\ref{T:variance1}.

\begin{proof}[Sketch of the proof of Theorem~\ref{T:variance1}]
Let $\phi: B_s(H) \to B_s(H)$ be a surjective linear isometry. Clearly,
$\phi$ preserves the extreme points of the unit ball of $B_s(H)$ which
are well-known (and easily seen) to be exactly the self-adjoint
unitaries, i.e.,
the operators of the form $2P-I$ where
$P$ is a projection.
Now, one can readily prove that among those extreme points, $I$
and $-I$ are distinguished by the following property. The extreme point
$U$ is either $I$ or $-I$ if and only if we have $\| U-V\| \in \{ 0,2\}$
for every extreme point $V$. Therefore, we get $\phi(\{ I,-I\})=\{
I,-I\}$. Clearly, we can suppose without loss of generality that
$\phi(I)=I$.
In that case we obtain that $\phi$ preserves the projections.
This gives us that $\phi$ is a Jordan automorphism of $B_s(H)$,
that is, it satisfies the equality
$\phi(AB+BA)=\phi(A)\phi(B)+\phi(B)\phi(A)$ for every $A,B\in B_s(H)$
(cf. \cite{BrSe4} or \cite{BrSe3}). Therefore,
we have that $\phi$ is of the form
\[
\phi(A)=UAU^* \qquad (A\in B_s(H))
\]
with some unitary or antiunitary operator $U$ on $H$ (see, for example,
\cite{CVLL}).
\end{proof}

The proof of Theorem~\ref{T:variance2} is much more difficult than
the one given above and is based on the following series of lemmas.
Our first observation below will prove to be fundamental from the
view-point of the proof of Theorem~\ref{T:variance2} that we are going
to present. It states that the
maximal deviation of an operator $T$ is equal to the so-called factor
norm of $T$ in the factor Banach space $B_s(H)/\mathbb R I$. (In
particular, this result implies that the function $T \mapsto \| T\|_v$
is a semi-norm on $B_s(H)$.)
Denote by ${\overline T}$ the equivalence class of
$T$ in $B_s(H)/\mathbb R I$. The factor norm $\| {\overline T}\|$ of
$T$ is defined by
\[
\| {\overline T}\| =\inf_{\lambda \in \mathbb R} \| T+\lambda I\|.
\]
As the spectral radius and the operator norm
of a self-adjoint operator are the same, it easily follows that
$\| {\overline T} \|$
is equal to the half of the diameter of the spectrum $\sigma(T)$ of $T$.

\begin{lemma}\label{L:var1}
For all $T\in B_s(H)$ we have $\|
T\|_v=\|{\overline T}\|=\text{diam}(\sigma(T))/2$.
\end{lemma}

\begin{proof}
As we have already verified that
$\|{\overline T}\|=\text{diam}(\sigma(T))/2$,
we have to prove only the first equality.
For a scalar operator $T$, both quantities $\| T\|_v$ and $\|
{\overline T}\|$ are 0.
Otherwise, we can assume that $ 0\leq T\leq I$ and that $\{0,1\}\subset
\sigma(T)\subset [0,1]$. This is because the factor norm and the
maximal deviation
of $T$ are invariant under adding scalar operators and they are absolute
homogeneous. In this case we have $\|{\overline T}\|=\frac{1}{2}$.

First we prove the easier inequality $\|T\|_v\leq\|{\overline T}\|$. For
any $\lambda\in\mathbb R$ we have
\begin{equation*}
\begin{gathered}
  \|T\|_v^2=\|T+\lambda I\|_v^2=\sup_{ \|\varphi\|=1}
  \left(\langle(T+\lambda I)^2\varphi,\varphi\rangle-\langle (T+\lambda
  I)\varphi,\varphi\rangle^2\right)
\leq\\
  \sup_{\|\varphi\|=1}\langle(T+\lambda I)^2\varphi,\varphi\rangle
  =\|(T+\lambda I)^2\|=\|T+\lambda I\|^2.
\end{gathered}
\end{equation*}
This yields $\|T\|_v\leq\|T+\lambda I\|$ for all $\lambda\in\mathbb R$
which implies that $\|T\|_v\leq\|{\overline T}\|$.

Now, we turn to the less obvious inequality $\frac{1}{2}=\|{\overline
T}\|\leq\|T\|_v$. Let $E_T$ be the spectral measure corresponding to
$T$. Since $0$ and $1$ are in the spectrum of $T$,
it follows that for any $0<\delta \leq
\frac{1}{2}$, the measures of $]-\delta,\delta[\cap\sigma(T)$
and
$]1-\delta,1+\delta [\cap\sigma(T)$ under $E_T$ are mutually orthogonal
nonzero projections.
At this stage $\delta$ is not fixed, we shall specify it later.
Denote these projections by $P_0$ and $P_1$, respectively.

Let $x$ be a unit vector in the range of $P_0$ and $y$ be a unit
vector in the range of $P_1$.
Define $\varphi =(x+y)/\sqrt{2}$. Then $\varphi\in H$ is a unit
vector and we assert that the following inequality holds
\begin{equation}\label{E:var1}
\sqrt{\langle T^2\varphi,\varphi\rangle-\langle T\varphi,\varphi\rangle^2}
\geq
\sqrt{\frac{(1-2\delta)^2}{2}-\frac{(1+2\delta)^2}{4}}.
\end{equation}
To see this, first observe that
$Tx=TP_0x$ and $Ty=TP_1y$.
Since
\[
TP_0=\int_{]-\delta, \delta[\cap \sigma(T)} t\,d\, E_T(t),
\]
we deduce $\| TP_0\| \leq \delta$. This yields that
\[
\| Tx\|\leq \delta .
\]
A similar argument shows that $\|Ty-y\|=\|TP_1y-P_1y\|\leq\delta$.
Since $\| y\|=1$, this gives us that
\[
1-\delta\leq\|Ty\|\leq 1+\delta .
\]
Now, to prove \eqref{E:var1} we estimate $\langle
T^2\varphi,\varphi\rangle= \|T\varphi\|^2$
from below and $\langle T\varphi,\varphi\rangle^2$ from above.
Since $T\varphi=(Tx+Ty)/\sqrt 2$, we have
\[
\|T\varphi\|
\geq
\frac{-\| Tx\|+\|Ty\|}{\sqrt 2}
\geq
\frac{-\delta+1-\delta}{\sqrt 2}
\]
and thus we get
\begin{equation}\label{E:var2}
\langle T^2\varphi,\varphi\rangle=
\|T\varphi\|^2\geq \frac{(1-2\delta)^2}{2}.
\end{equation}
Using the equality $TP_0=P_0T$ and the fact that
$P_0$ and $P_1$ are mutually orthogonal projections, we have
\[
\langle Tx,y\rangle=\langle TP_0x,P_1y\rangle=
\langle P_0Tx,P_1y\rangle=\langle Tx,P_0P_1y\rangle=0.
\]
This also implies that $\langle Ty,x\rangle=0$.
Therefore, we infer
\[
\langle T\varphi,\varphi\rangle=
\frac{1}{2}\left(\langle Tx,x\rangle+\langle Ty,y\rangle\right).
\]
Since $|\langle Tx,x\rangle|\leq\|Tx\|\leq\delta$ and
$|\langle Ty,y\rangle|\leq\|Ty\|\leq 1+\delta$, we obtain
\[
\langle T\varphi,\varphi\rangle^2\leq \frac{(1+2\delta)^2}{4}.
\]
This inequality together with \eqref{E:var2} gives \eqref{E:var1}.

Now, for an arbitrary $\epsilon >0$, choosing $\delta$ such that it
satisfies
\[
\sqrt{\frac{(1-2\delta)^2}{2}-\frac{(1+2\delta)^2}{4}}
\geq\frac{1}{2}-\epsilon,
\]
it follows from what we have already proved that
we can pick a unit vector $\varphi \in H$ for which
\[
\| T\|_v\geq
\sqrt{\langle T^2\varphi,\varphi\rangle-\langle T\varphi,\varphi\rangle^2}
\geq
\frac{1}{2}-\epsilon.
\]
This gives us that $\| T\|_v\geq \frac{1}{2}=\| {\overline T} \|$ which
completes the proof of the lemma.
\end{proof}

\begin{remark}\label{R:var1}
As we have seen,
the quantity $\| T\|_v=\| {\overline T}\|$ is exactly
the half of the diameter of the spectrum of $T$.
Therefore, if $T\geq 0$ and $0 \in \sigma(T)$, then
$\|T\|_v=\|{\overline T}\|\leq\frac{1}{2}$ if and only if $0\leq T\leq
I$.
\end{remark}

This observation will be used in the proof
of our next lemma which determines the extreme points of the
(closed) $\frac{1}{2}$-ball of the Banach space $B_s(H)/\mathbb R I$.

\begin{lemma}\label{L:var2}
The extreme points of the ball
$\{{\overline A}\in B_s(H)/\mathbb R I \,: \, \|{\overline
A}\|\leq\frac{1}{2}\}$
are the classes of nontrivial projections, that is,
the elements ${\overline P}\in B_s(H)/\mathbb R I$, where $P$ is a
nontrivial projection ($P\neq 0,I$) on $H$.
\end{lemma}

\begin{proof}
The point in the proof is to reduce the problem concerning classes of
operators to a problem concerning single operators.

First we check that the classes of nontrivial projections
are extreme points of the ball in question.
Suppose that $P$ is a nontrivial projection and
\[
{\overline P}=\mu{\overline T}+(1-\mu){\overline S},
\]
where $0<\mu<1$, $\| {\overline T}\|\leq\frac{1}{2}$,
$\| {\overline S}\|\leq\frac{1}{2}$, $T,S\in B_s(H)$.
Adding scalar operators if necessary,
we can suppose that $T,S\geq 0$, $0\in \sigma(T)$, $0\in \sigma(S)$.
Clearly,
\[
P=\mu T+(1-\mu)S+\lambda I
\]
holds for some $\lambda \in \mathbb R$.

We claim that $\lambda =0$.
If $\varphi\in H$ is a unit vector in the kernel of $P$, we infer that
\[
0=\langle P\varphi,\varphi\rangle=
\mu \langle T\varphi,\varphi\rangle+
(1-\mu) \langle S\varphi,\varphi\rangle+\lambda.
\]
Since $\langle T\varphi,\varphi\rangle\geq 0$ and $\langle
S\varphi,\varphi\rangle\geq 0$, the above equality yields $\lambda\leq
0$.

It follows from $\sigma(P)=\{0,1\}$ that
$\|{\overline P}\|=\frac{1}{2}$.
We compute
\begin{equation*}
  \frac{1}{2}=\|{\overline P}\|=\|\mu
{\overline T}+(1-\mu){\overline S}\|
  \leq
\mu\|{\overline T}
\|+(1-\mu)\|{\overline S}\|\leq(\mu+1-\mu)\frac{1}{2}=\frac{1}{2},
\end{equation*}
from which we deduce that $\|{\overline
T}\|=\|{\overline S}\|=\frac{1}{2}$.
Using Remark~\ref{R:var1} we get $0\leq T,S\leq I$.
So, if $\varphi$ is a unit vector in the range of
$P$, then we have
\[
1=\langle P\varphi,\varphi\rangle=
\mu\langle T\varphi,\varphi\rangle+(1-\mu)\langle S\varphi,\varphi\rangle+
\lambda\leq \mu+(1-\mu)+\lambda,
\]
which gives us that $\lambda\geq 0.$ Therefore, it follows that
$\lambda=0$ as we have claimed.

Consequently, we have $P=\mu T+ (1-\mu)S$.
This means that $P$ is a nontrivial convex combination of two elements
of the operator interval $[0,I]$. However, it is well-known that the
extreme points of this operator interval
are exactly the projections. Hence, we get $P=T=S$.
This proves that the classes of nontrivial projections are
really extreme points.

It remains to prove that these classes are the only extreme points.
In order to see this, let
$B$ be a self-adjoint operator with $\| {\overline B}\|= \frac{1}{2}$
which is not a nontrivial projection.
We show that ${\overline B}$ is not an extreme point of the ball in question.
Clearly, just as above, we can assume that $B\geq 0$ and
$0\in\sigma(B)$.
Then we have $0\leq B\leq I$.
As $\| {\overline B}\|=\frac{1}{2}$, it also follows that $1\in
\sigma(B)$.
We are going to show that there exist two operators $B_1, B_2$
in the operator interval $[0,I]$ such that $B=(B_1+B_2)/2$ and $B\neq
B_1,B_2$.
In the present situation this will imply that ${\overline B} \neq
\overline{B_1}, \overline{B_2}$. Then, as
${\overline{B}}=({\overline{B_1}}+{\overline{B_2}})/2$,
$\|{\overline{B_1}}\|,\|{\overline{B_2}}\| \leq \frac{1}{2}$ (see
Lemma~\ref{L:var1}), we can infer that ${\overline B}$
is not an extreme point. So, in order to construct such operators
$B_1,B_2$, choose $\lambda_0\in \sigma(B)\cap ]0,1[$. (The existence of
such a $\lambda_0$ follows from the facts that $B$ is not
a non-trivial projection and that $\|{\overline B}\|\neq 0$.)
Now, one can easily find continuous real valued functions
$f_1,f_2:[0,1]\to
[0,1]$ such that $(f_1+f_2)/2$ is the identity on $[0,1]$ and
$f_1(\lambda_0)\neq \lambda_0 \neq f_2(\lambda_0)$. Defining
$B_1=f_1(B), B_2=f_2(B)$,
it follows from the properties of the continuous function calculus
that we obtain operators with the desired properties.
This completes the proof of the lemma.
\end{proof}

In what follows, we intend to characterize the unitary equivalence of
nontrivial projections $P,Q$ by means of some correspondence between
the classes ${\overline P}$ and ${\overline Q}$ that can be expressed
in terms of the metric induced by the factor norm. The first step in
this direction is made in the following lemma.

\begin{lemma}\label{L:var3}
Let $P$ and $Q$ be projections on $H$. Suppose that $P$ is nontrivial
and $\|{\overline P}-{\overline Q}\|<\frac{1}{2}$. Then $P$ is unitarily equivalent to
$Q$.
\end{lemma}

\begin{proof}
First observe that $Q\neq 0,I$.
In fact, in the opposite case we would have $\| {\overline P}\|<1/2$.
But this means that the diameter of $\sigma(P)$ is less
than 1, which gives us that $P$ is a trivial projection, a
contradiction.

Because of the definition of the factor norm there exists a
$\mu\in\mathbb R$ such that $\|P-(Q+\mu I)\|<\frac{1}{2}$.
Let $R$ be a projection of rank at most 2 whose
range contains a unit vector from the range of $P$ and a unit vector
from the range of $Q$, respectively.
The operators $RPR$ and $R(Q+\mu I)R$ are of finite rank,
$1$ is the largest eigenvalue of $RPR$ and $1+\mu$ is
the largest eigenvalue of $R(Q+\mu I)R$.
Indeed, to prove for example this last statement, observe that
\[
R(Q+\mu I)R\leq R(I+\mu I)R=(1+\mu)R\leq (1+\mu)I.
\]
This shows that
the spectrum of $R(Q+\mu I)R$ is a subset of the interval $]-\infty,
1+\mu ]$. On the other hand, $1+\mu$ is an eigenvalue of the operator
$R(Q+\mu I)R$
since the range of $R$ contains a unit vector from the range of $Q$.

By Weyl's perturbation theorem
(see, for example, \cite[Corollary III.2.6]{Bha})
we deduce
\[
\begin{gathered}
|\mu|=|1-(1+\mu)|\leq\|RPR-R(Q+\mu I)R\| \\
\leq \|R\|\|P-(Q+\mu I)\|\|R\|<\frac{1}{2},
\end{gathered}
\]
and so we have
\[
\|P-Q\|\leq\|P-(Q+\mu I)\|+|\mu|<\frac{1}{2}+\frac{1}{2}=1.
\]
But it is a well-known result that if the distance between two
projections
in the operator norm is less than 1, then they are unitarily equivalent.
This completes the proof of the lemma.
\end{proof}

A useful solution of the problem mentioned before Lemma~\ref{L:var3} is
given in the next result.

\begin{lemma}\label{L:var4}
Let $P$ and $Q$ be projections on $H$ and suppose that $P$ is
nontrivial.
Then $P$ is unitarily equivalent to $Q$ if and only if there
exists a continuous function $\varphi :[0,1]\to \overline{P(H)}$
such that $\varphi(0)={\overline P}$ and $\varphi(1)={\overline Q}$.
\end{lemma}

(Here $\overline{P(H)}$ denotes the set of classes in $B_s(H)/\mathbb R
I$ which correspond to projections.)

\begin{proof}
The necessity is easy to see. Indeed, this follows from the
well-known fact that if $P,Q$ are equivalent projections then they can
be connected by a continuous curve (continuity is meant in the operator
norm topology) in the set of projections and from the fact that the
operator norm majorizes the factor norm.

Now, conversely, suppose that there exists a continuous mapping
$\varphi:[0,1]\to\overline{P(H)}$ such that $\varphi(0)={\overline P}$
and $\varphi(1)={\overline Q}$. As $\varphi$ is defined on a compact set,
it is
uniformly continuous. Hence, we can choose a positive $\delta$
such that
\[
\|\varphi(t)-\varphi(s)\|<\frac{1}{2}
\qquad\text{if}\quad
|s-t|<\delta,\,\,\,s,t\in[0,1].
\]
It follows that there exist projections $P_1,\dots,P_n$ with the
property that
\[
\|{\overline P}-\overline{P_1}\|<\frac{1}{2}
\quad , \ldots , \quad
\|\overline{P_n}-{\overline Q}\|<\frac{1}{2}.
\]
By Lemma~\ref{L:var3}, we obtain that $P$ and $P_1$ are
unitarily equivalent (and, consequently, $P_1$ is non-trivial). Using
this argument again and again we can conclude that $P$
is unitarily equivalent to $Q$.
\end{proof}

The meaning of our last lemma which follows is a metric characterization
of the equality of nontrivial projections in $B_s(H)$ with respect to
the semi-norm $\| .\|_v$. Denote by
$F_s(H)$ the set of all finite rank elements in $B_s(H)$.

\begin{lemma}\label{L:var5}
Let $P$ and $Q$ be nontrivial projections on $H$ such that
\[
\|P+A\|_v=\|Q+A\|_v
\]
holds for all $A\in F_s(H)$. Then we have $P=Q$.
\end{lemma}

\begin{proof}\footnote{We remark that in his/her report the
referee presented a more elementary proof of this lemma which
uses only matrix (finite dimensional) arguments.} Let $R$ be a
rank-1 subprojection of the projection $P$. Then the diameter of
the spectrum of $P+R$ is 2, so by Lemma~\ref{L:var1} we have
\[
1=\|P+R\|_v=\|Q+R\|_v,
\]
that is, the diameter of $\sigma(Q+R)$ is
also equal to $2$. Since $0\leq Q+R\leq 2I$, thus $\sigma(Q+R)$ is
a subset of the closed interval $[0,2]$. Therefore, we have $0,2\in
\sigma(Q+R)$.

It is well-known that the spectrum of any normal operator
coincides with its approximate point spectrum. Consequently, we can find
unit vectors $x_n$ in $H$ $(n\in \mathbb N)$ such that
\[
\|Qx_n+Rx_n-2x_n\|\to 0\quad\text{as}\quad n\to\infty.
\]
This yields that
\begin{equation}\label{E:var3}
\|Qx_n+Rx_n\|\to 2.
\end{equation}
Denote $u_n=Qx_n$ and $v_n=Rx_n$.
We have $\|u_n\|\leq 1,\|v_n\|\leq 1$. Since $v_n$ is in
the range of $R$ which is 1-dimensional, there must exist a convergent
subsequence of $(v_n)$.
Without any loss of generality we can assume that
this subsequence is $(v_n)$ itself. So, there
exists a vector $v$ in the range of $R$ such that $\|v_n-v\|\to 0$.
Since
\[
\big| \|u_n+v\|-\|u_n+v_n\| \big|
\leq
\|v-v_n\|\to 0
\]
and $\| u_n+v_n\|\to 2$, we have $\|u_n+v\|\to 2$.
On the other hand, by the parallelogram identity we obtain
\[
\|u_n-v\|^2=2\|u_n\|^2+2\|v\|^2-\|u_n+v\|^2.
\]
Therefore, we have
\[
\limsup_{n\to\infty}\|u_n-v\|^2\leq 2+2-4=0,
\]
which implies that $\|u_n -v\| \to 0$.
So, both $(u_n),(v_n)$ converge to $v$. Taking \eqref{E:var3} into
account, it is clear that $v\neq 0$.

Since the sequence $(u_n)$ is in the range of $Q$ which
is a closed subspace, it follows that its limit $v$ also belongs to this
range.
But $v$ generates the range of $R$ and hence $R$ is a subprojection of
$Q$. So, we have proved the
following: every rank-1 subprojection of $P$ is a subprojection of $Q$.
Therefore, $P$ is a subprojection of $Q$.
Changing the role of $P$ and $Q$, we get that $Q$ is also a
subprojection of $P$ and hence we obtain $P=Q$.
\end{proof}

Now, we are in a position to prove our main result.

\begin{proof}[Proof of Theorem~\ref{T:variance2}]
The brief summary of the proof is as follows.
Our transformation $\phi$ which preserves the maximal deviation induces
a surjective linear isometry $\Phi$ on the factor space $B_s(H)/\mathbb
R I$.
This $\Phi$ necessarily preserves the extreme points of the
$\frac{1}{2}$-ball
which points are well characterized in Lemma~\ref{L:var2}. This implies
a certain preserving property of the original transformation $\phi$.
Namely, we obtain that $\phi$ preserves the operators of the form
"nontrivial projection + scalar$\cdot I$". This will imply
that $\phi$ preserves the commutativity on $F_s(H)+\mathbb R I$.
Extending $\phi$ from this set to its complex linear span $F(H)+\mathbb
C I$ ($F(H)$ stands for the set of all finite rank bounded linear
operators on $H$), we obtain a
complex-linear transformation which preserves normal operators. Applying
the technique of the proof of a nice result of Bre\v sar and \v Semrl
given in
\cite{BrSe1}, we can conclude the proof in the case when $\dim H\geq 3$.
If $\dim H=2$, then rather surprisingly we can reduce our problem
quite easily to Wigner's classical unitary-antiunitary theorem.
So, this is the plan what we now carry out.

Define a map $\Phi:B_s(H)/\mathbb R I\to B_s(H)/\mathbb R I$
in the following way
\[
\Phi({\overline A})=\overline{\phi(A)} \qquad (A\in B_s(H)).
\]
The transformation $\phi$ is a linear bijection of $B_s(H)$ which
preserves the maximal deviation.
By Lemma~\ref{L:var1}, we easily obtain that $\phi$ preserves the scalar
operators and then that $\Phi$ is a well-defined
linear bijection on $B_s(H)/\mathbb R I$ which preserves the factor
norm. It follows that $\Phi$ preserves all closed balls around
${\overline 0}$
as well as their extreme points. Therefore, by Lemma~\ref{L:var2}, we
deduce that $\phi$ preserves the set of all operators of the form
$P+\lambda I$, where $P$ is a nontrivial projection and $\lambda\in
\mathbb R$.

We shall show that $\phi$ preserves the commutativity on $F_s(H)+\mathbb
R
I$. Let $P'$ and $Q'$ be mutually orthogonal projections. We known that
there exist projections $P,Q,R$ and real numbers
$\lambda_1,\lambda_2,\lambda_3$ such that
\begin{align*}
 \phi(P')&=P+\lambda_1 I\\
 \phi(Q')&=Q+\lambda_2 I\\
 \phi(P'+Q')&=R+\lambda_3 I.
\end{align*}
By the linearity of $\phi$
this implies that $P+Q=R+tI$ for some real number $t$ (in fact,
$t=\lambda_3-\lambda_1-\lambda_2$).
We assert that $P$ and $Q$ are either commuting or the
projections $P,Q, R$ are unitarily equivalent to each other.

In order to prove this, we distinguish the following cases.

{\sc Case I.} Suppose that $R$ is scalar. Then
$P+Q$ is also scalar which implies that $P,Q$ commute.

{\sc Case II.} Suppose that $R$ is not scalar, that
is, $R$ is a nontrivial projection. Consider the orthogonal
decomposition of $H$ induced by the range and the kernel of $R$.
Every operator has a matrix representation with respect to this
decomposition. As for $P+Q$, we can write
\begin{equation}\label{E:var21}
P+Q=R+tI=\begin{bmatrix}
         (1+t)I & 0\\
         0 & tI
         \end{bmatrix}.
\end{equation}
The inequality $0\leq P+Q\leq 2I$ implies that $0\leq t\leq 1$.
According to the possible values of $t$ we have the following
sub-cases.

{\sc Case II/1.} Suppose that $t=0$. Then $P+Q=R$ is a projection
and hence $(P+Q)^2=P+Q$. From this equality we easily deduce $PQ=QP=0$
which implies that $P,Q$ commute.

{\sc Case II/2.} Suppose that $t=1$. Then $P+Q=R+I$, which
implies that $R+(I-Q)=P$ is a projection. Just as above, we obtain that
$R, I-Q$ are commuting projections.
This implies that $R,Q$ commute and, finally, it follows from the
equality $R+(I-Q)=P$ that $P,Q$ also commute.

{\sc Case II/3.}\footnote{The main part of the argument in this
case was suggested by the referee. Due to his/her idea, the
original proof could be reduced from 3 pages to few lines.}
Suppose that $0< t<1$. In this case we use the result that any
two projections in generic position (i.e., with no common
eigenvectors) are unitarily equivalent (see \cite{Dav, Hal}). As
the spectrum of $P+Q=R+tI$ is contained in $\{t, 1+t\}$, the
numbers 0,1,2 are not in the spectrum of $P+Q$. This implies that
$P,Q$ are in a generic position and hence they are unitarily
equivalent. Similarly, as the spectrum of $R-P$ is contained in
$\{ -t, 1-t\}$ which does not contain -1,0,1, we infer that $P,R$
are in a generic position and hence they are unitarily
equivalent. It follows that the projections $P,Q,R$ are pairwise
unitarily equivalent. What does this mean for our original
projections $P',Q'$? Obviously, in the present case $P,Q,R$ are
nontrivial. Using Lemma~\ref{L:var4} and the isometric property
of $\Phi$ with respect to the factor norm, we obtain that the
projections $P',Q',P'+Q'$ are pairwise unitarily equivalent. But
if $P',Q'$ are nonzero mutually orthogonal finite rank
projections, then this can not happen.

Therefore, we have proved that for any finite rank projections $P',Q'$
with $P'Q'=Q'P'=0$ it follows that $\phi(P')\phi(Q')=\phi(Q')\phi(P')$.
If we pick operators $A,B\in F_s(H)$ which commute, then they
can be
diagonalized simultaneously. Using the just proved property of $\phi$
one can easily deduce that $\phi(A), \phi(B)$ also commute.

We show that
\[
\phi(F_s(H)+\mathbb R I)=
F_s(H)+\mathbb R I.
\]
If $\dim H<\infty$, this is obvious. So, let $H$
be infinite dimensional.
Pick a nonzero finite rank projection $P'$. Then
$\phi(P')=P+\lambda I$ holds for some nontrivial projection $P$ and
real number $\lambda$. If $P$ is of finite rank or of finite corank,
then we obtain $\phi(P')\in F_s(H)+\mathbb R I$. So, let us see
what happens if $P$ is of infinite rank and infinite corank.

First suppose that $\dim \text{rng}\, P\leq \dim \text{rng}\, P^\perp$.
Then we can find nontrivial projections $P_1$ and $P_2$ such that
$P=P_1+P_2$ and $P,P_1,P_2$ are mutually unitarily equivalent. Now,
referring to Lemma~\ref{L:var4}, there are nontrivial projections
$P_1',P_2'$ such that
\[
P'+\mu I=P_1'+P_2'
\]
holds for some $\mu \in \mathbb R$ and the projections $P',P_1',P_2'$
are mutually unitarily
equivalent. So, the projections $P_1',P_2'$ are of finite rank  and we
see that on the right hand side of the equality above there is a
finite rank operator. This gives us that $\mu$ must be zero and then we
have $P'=P_1'+P_2'$. Like in the argument given in {\sc Case II/1.}, we
obtain
that $P_1',P_2'$ are mutually orthogonal projections. We now conclude
that, because of unitary equivalence and orthogonality,
the equality $P'=P_1'+P_2'$ is untenable which is a contradiction.

Next suppose that
$\dim \text{rng}\, P\geq \dim \text{rng}\, P^\perp$.
Then we can apply the argument above for $P^\perp$ to
find nontrivial projections $P_1$ and $P_2$ such
that $P^\perp=P_1+P_2$ and $P^\perp,P_1,P_2$ are mutually unitarily
equivalent. This implies that there are
nontrivial projections $P_1',P_2'$ such that
\begin{equation}\label{E:var31}
{P'}^\perp+\nu I=P_1'+P_2'
\end{equation}
holds for some $\nu \in \mathbb R$ and the
projections ${P'}^\perp,P_1',P_2'$ are mutually unitarily equivalent.
(Observe that, as
$\Phi({\overline {P'}})={\overline P}$, we have
$\Phi({\overline{{P'}^\perp}})={\overline{P^\perp}}$.)
It follows that
the projections $P_1',P_2'$ are of finite corank and hence their ranges
have nonempty intersection. Therefore, we obtain that 2 belongs to the
spectrum of the operator $P_1'+P_2'$, and by \eqref{E:var31} this
implies that $\nu=1$. Now, the equation \eqref{E:var31} can be rewritten
in the form
\[
P'=(I-P_1')+(I-P_2')={P_1'}^\perp+{P_2'}^\perp,
\]
where the nontrivial projections
$P',{P_1'}^\perp,{P_2'}^\perp$ are pairwise unitarily equivalent.
Just as in the previous paragraph we arrive at a contradiction.

Therefore, we have $\phi(P')\in F_s(H)+\mathbb R I$ for every finite
rank projection $P'$. Applying the spectral theorem for self-adjoint
finite rank operators,
it follows that $\phi(F_s(H)+\mathbb R I)\subset F_s(H)+\mathbb R I$.
As $\phi^{-1}$ has the same properties as $\phi$, considering the above
relation for $\phi^{-1}$ in the place of $\phi$, we conclude that
\[
\phi(F_s(H)+\mathbb R I)= F_s(H)+\mathbb R I.
\]

To sum up what we have already proved, it has turned out that $\phi$
when restricted onto $F_s(H)+\mathbb R I$ is a bijective linear map
which preserves commutativity.
Consider the complex unital algebra $F(H)+\mathbb C I$. As the real
and imaginary parts of an operator in $F(H)+\mathbb C I$ belong to
$F_s(H)+\mathbb R I$, one can readily verify that the map $\tilde
\phi: F(H)+\mathbb C I\to F(H)+\mathbb C I$ defined by
\[
\tilde \phi(A+iB)=\phi(A)+i\phi(B) \qquad (A,B\in F_s(H)+\mathbb R I)
\]
is a bijective complex-linear
transformation. It is an elementary fact that a bounded linear operator
is normal if and only if its real and imaginary parts are commuting.
As $\phi$ preserves commutativity between self-adjoint finite rank
operators, it follows that $\tilde \phi$ preserves
normality.
If $\dim H\geq 3$, then this latter preserving property is strong enough
to imply that $\tilde \phi$ is of a certain particular form. In fact,
there is a nice result of
Bre\v{s}ar and \v{S}emrl \cite[Theorem 2]{BrSe1}
which, in the case when $\dim H\geq 3$, characterizes the bijective
linear mappings on $B(H)$ that preserve normal operators.
Although the algebra on which our transformation $\tilde \phi$ is
defined differs from $B(H)$ in general, it is not hard to see that
the technique used in \cite{BrSe1} can be applied to our
present situation as well.
This gives us the following two possibilities for the form of $\tilde
\phi$:
\begin{enumerate}
\item[(i)]  there exist a unitary operator $U$ on $H$, a linear
            functional $f:F(H)+\mathbb C I\to\mathbb C$ and
            a scalar $c\in\mathbb C$ such that
            \[
               \tilde \phi(T)=c UTU^*+f(T)I \qquad (T \in
               F(H)+\mathbb C I)
            \]
\item[(ii)] there exist an antiunitary operator $U$ on $H$,
            a linear functional $f:F(H)+\mathbb C I\to\mathbb C$ and
            a scalar $c\in\mathbb C$ such that
            \[
               \tilde \phi(T)=c UT^*U^*+f(T)I \qquad (T \in
               F(H)+\mathbb C I).
            \]
\end{enumerate}
Concerning $\phi$, this means that there is an either unitary or
antiunitary
operator $U$ on $H$, a real-linear function $f:F_s(H)+\mathbb R I\to
\mathbb C$, and a constant $c\in \mathbb C$ such that
\[
\phi(A)=c UAU^*+f(A)I \qquad (A \in F_s(H)+\mathbb R I).
\]
As $\phi(A)$ is self-adjoint, we have
\begin{equation}\label{E:var32}
\overline{c} UAU^*+\overline{f(A)}I=c UAU^*+f(A)I
\end{equation}
for every $A\in F_s(H)+\mathbb R I$. If $A$ is not a scalar operator,
then it follows from this
equality that $\overline{c}=c$. Next, we obtain from \eqref{E:var32}
that
$f$ is real valued. As $\phi$ preserves maximal deviation, we obtain
that $|c|=1$. Therefore, $c=\pm 1$ and we have the desired form for our
transformation $\phi$ on $F_s(H)+\mathbb R I$.
It remains to show that the same formula holds also on the whole
space $B_s(H)$.

In order to see this, observe that
composing $\phi$ by the transformation $A\mapsto cU^* AU$, we can assume
without loss of generality that
\[
\phi(A)=A+l(A) I
\]
holds for every $A\in F_s(H)+\mathbb R I$,
where $l:F_s(H)+\mathbb R I\to \mathbb R$ is a linear functional.
Let $P$ be a nontrivial projection on $H$. We know that $\phi(P)=Q+\mu
I$ for some nontrivial projection $Q$ and real number $\mu$.
Pick an arbitrary $A\in F_s(H)$. Since $\phi(A)$ is a scalar
perturbation of $A$, we have
\[
\| Q+A\|_v = \|\phi(P)+A\|_v =
\|\phi(P)+\phi(A)\|_v = \|\phi(P+A)\|_v = \|P+A\|_v.
\]
Since this holds true for every self-adjoint finite rank operator $A$,
it follows from Lemma~\ref{L:var5} that $Q=P$. This gives us that
$\phi(P)-P\in \mathbb R I$ which holds also in the case when $P$ is
trivial.
So, we have $\Phi({\overline P})={\overline P}$ for every projection
$P$. Since the linear transformations $A\mapsto \Phi({\overline A})$ and
$A\mapsto {\overline A}$ are continuous (on $B_s(H)$ we consider
the operator norm while $B_s(H)/\mathbb R I$ is equipped with the factor
norm), they are equal on the projections, it follows from the
spectral theorem of self-adjoint operators and from the properties of
the spectral integral that we have
$\Phi({\overline A})={\overline A}$ for every $A\in B_s(H)$.
This gives us that
\[
\phi(A)-A \in \mathbb R I \qquad (A\in B_s(H))
\]
which obviously implies that there is a linear functional $h:B_s(H)
\to \mathbb R$ such that
\[
\phi(A)=A +h(A) I \qquad (A\in B_s(H)).
\]
This completes the proof in the case when $\dim H\geq 3$.

As the statement of the theorem is trivial for $\dim H=1$, it remains to
consider
the case when $\dim H=2$. In this case the nontrivial projections
are exactly the rank-one projections. Pick a rank-one projection $P$. We
know that there is a rank-one projection $P'$ such that $\phi(P)$ is
equal to the sum of $P'$ and a scalar operator. It is easy to see that
this $P'$ is unique.
(In fact, one can prove independently from the dimension of $H$
that in the class of every nontrivial projection there is only
one projection.)
Therefore, we can denote $P'=\psi(P)$ and
obtain a bijective transformation $\psi$ on the set of all rank-one
projections. We assert that $\psi$ has the property that
\begin{equation}\label{E:var5}
\text{tr}\, PQ=\text{tr}\, \psi(P)\psi(Q)
\end{equation}
holds for arbitrary rank-one projections $P,Q$ on $H$. Here $\text{tr}$
denotes the usual trace functional. As $\phi$ preserves
the maximal deviation, this will clearly follow from the equality
\begin{equation}\label{E:var33}
\| P-Q\|_v =\sqrt{1-\text{tr}\, PQ}
\end{equation}
that we are going to prove now.
In fact, observe that the maximal
deviation and the trace functional are invariant under the
transformations $A\mapsto VAV^*$, where $V$ is any unitary operator.
Therefore, we can assume that
\[ P=
       \begin{bmatrix}
             1 & 0\\
             0 & 0
        \end{bmatrix}
\]
while $Q$ is an arbitrary self-adjoint idempotent 2 by 2 matrix.
It is easy to check that $Q$ is of the form
\[
Q=
       \begin{bmatrix}
             a & \sqrt{a(1-a)}e^{i\theta}\\
             \sqrt{a(1-a)}e^{-i\theta} & 1-a
        \end{bmatrix}
\]
where $a, \theta$ are real numbers and $0\leq a\leq 1$.
We have that the eigenvalues of $P-Q$ are $\pm
\sqrt{1-a}$ and hence obtain that
$\| P-Q\|_v=\sqrt{1-a}$. On the other hand, it is trivial to check that
$\text{tr}\, PQ=a$. This results in the desired equality
\eqref{E:var33}.

So, we have a bijective transformation $\psi$ on the set of all
rank-one projections which satisfies \eqref{E:var5}. Wigner's
classical theorem on quantum mechanical symmetries (the so-called
unitary-antiunitary theorem) describes the form of exactly such
transformations in the case of general Hilbert spaces. We obtain
that there exists an either unitary or antiunitary operator $U$
on $H$ such that
\[
\psi(P)=UPU^*
\]
holds for every rank-one projection $P$. As $\phi(P)$ differs
from $\psi(P)$ only by a scalar operator, we obtain that
$\phi(P)-UPU^*\in \mathbb R I$. By
linearity this gives us that $\phi(A)-UAU^*$ is a scalar operator for
every $A\in
B_s(H)$. Now, one can easily complete the proof
in the case when $\dim H=2$.
\end{proof}

\begin{remark}
As it is seen, preserving commutativity has played important role in our
proof above. In fact, preserver problems of this kind are among the most
fundamental ones in the theory of LPP's. To mention one of the most
well-known results of this type which concerns operator algebras, we
refer to \cite{Oml}
\end{remark}

\begin{proof}[Proof of Theorem~\ref{T:variance3}]
This follows immediately from Theorem~\ref{T:variance1}
using the following important result of Mazur and Ulam \cite{MaUl}.
If $\mathcal V$ is a real normed vector space and $T:\mathcal V\to
\mathcal V$ is a bijective map which preserves the distance on
$\mathcal V$
(i.e., $T$ satisfies $\| T(x)-T(y)\|=\| x-y\|$ $(x,y\in \mathcal V)$),
then
$T$ can be written in the form $T(x)=L(x)+x_0$ $(x\in \mathcal V)$,
where $L:\mathcal V\to
\mathcal V$ is a bijective linear isometry and $x_0\in \mathcal V$ is a
fixed vector.
\end{proof}

As for the proof of Theorem~\ref{T:variance4}, we have to work
more than in the previous proof as $\| .\|_v$ is only a semi-norm.

\begin{proof}[Proof of Theorem~\ref{T:variance4}]
Considering the map $A\mapsto \phi(A)-\phi(0)$, it is obvious that we
can assume that $\phi$ sends 0 to 0. In what follows we use this
assumption.

Consider the linear functional $\lambda I\mapsto \lambda$ on $\mathbb R
I$. Extend it to a linear functional $l$ of the whole vector space
$B_s(H)$. (We do not need any kind of continuity of $l$,
so no need to use Hahn-Banach theorem.)
Define the transformation $\phi_1:B_s(H)\to B_s(H)$ in the
following way
\[
\phi_1(A)=\phi(A)-l(\phi(A))I+l(A)I \qquad (A\in B_s(H)).
\]
We assert that $\phi_1:B_s(H)\to B_s(H)$ is a bijective linear map, it
preserves the distance (with respect to the semi-metric $d_v$) and for
every $A\in B_s(H)$, $\phi(A)$
and $\phi_1(A)$ differs only in a scalar operator. If this is really the
case, then we can apply Theorem~\ref{T:variance2} for $\phi_1$ and we
are done. So, it remains to
prove that $\phi_1$ has the mentioned properties. As the last two ones
are obvious from the definition, we have to prove only that $\phi_1$ is
linear and bijective.
We begin with the linearity.
As $\phi$ preserves the distance with respect to $d_v$ and we have
supposed that $\phi(0)=0$, it follows that $\phi$ preserves
the scalar operators (in fact, scalar operators
can be characterized by the equality $\| A\|_v=0$; see
Lemma~\ref{L:var1}). Next, it is easy to show that the formula
\[
\Phi(\overline{A})={\overline{\phi(A)}} \qquad (A\in B_s(H))
\]
defines a bijective isometry (distance preserving map) on
$B_s(H)/\mathbb R I$ with respect to the factor norm. We only prove the
isometric property. Indeed,
\[
\| \Phi({\overline A})-\Phi({\overline B})\|=
\| \phi(A)-\phi(B)\|_v=
\| A-B\|_v=
\| {\overline A}-{\overline B}\|
\]
holds for every $A,B\in B_s(H)$. Since $\Phi({\overline
0})={\overline{\phi(0)}}={\overline 0}$, by Mazur-Ulam theorem we obtain
that $\Phi$ is linear.
Thus, for any $A,B\in B_s(H)$ we have
\[
\Phi({\overline A}+{\overline B})=\Phi({\overline A})+\Phi({\overline
B}),
\]
that is,
\[
{\overline{\phi(A+B}})={\overline{\phi(A)+\phi(B)}}.
\]
This gives us that
$\phi(A+B)- (\phi(A)+\phi(B))$ is a scalar operator, say
\[
\phi(A+B)- (\phi(A)+\phi(B))=\lambda I.
\]
We compute
\[
\begin{gathered}
\phi(A+B)- (\phi(A)+\phi(B))=\lambda I\\
=l(\lambda I)I= l(\phi(A+B)- (\phi(A)+\phi(B)))I.
\end{gathered}
\]
This implies that
\[
\phi(A+B)-l(\phi(A+B))I=
\phi(A)-l(\phi(A))I+ \phi(B)-l(\phi(B))I.
\]
Adding $l(A+B)I=l(A)I+l(B)I$ to this equality,
we obtain the additivity of $\phi_1$. The
homogeneity can be proved in a similar way.

We next show that $\phi_1$ is injective. Suppose that
\[
0=\phi_1(A)=
\phi(A)-l(\phi(A))I+l(A)I
\]
holds for some $A\in B_s(H)$. Then $\phi(A)$ is a scalar operator, say
$\phi(A)=\lambda I$, and this implies that $A$ is also scalar, say
$A=\mu I$. It follows from the above equation that
\[
0=\lambda I-l(\lambda I)I+l(\mu I)I=(\lambda -\lambda +\mu)I
\]
which yields $\mu=0$, i.e., we have $A=0$. This proves the injectivity
of $\phi_1$.

Finally, we prove that $\phi_1$ is surjective. To show this, first
observe that, by the definition of $\phi_1$ and the surjectivity of
$\phi$, the range of $\phi_1$ and $\mathbb R I$
generate the whole space $B_s(H)$. So, if $\phi_1$ is not surjective,
then
we have $\text{rng}\, \phi_1\cap \mathbb R I=\{ 0\}$. This
means that the only scalar operator in the range of $\phi_1$ is 0.
Now, as $\phi(I)$ is a scalar operator, it follows that $\phi_1(I)$ is
also scalar.
As $\phi_1(I)\in \text{rng}\, \phi_1$, we obtain that $\phi_1(I)=0$,
which, by the injectivity of $\phi_1$ implies that $I=0$, a
contradiction. Therefore, $\phi_1$ must be surjective.
So, we have proved all the asserted
properties of $\phi_1$ and hence the proof of the theorem is complete.
\end{proof}

\section{An Open Problem}

To conclude the paper we give another interpretation of our main result
Theorem~\ref{T:variance2}. Namely, in view of Lemma~\ref{L:var1}, our
theorem
describes the form of all bijective linear transformations of $B_s(H)$
which preserve the diameter of the spectrum. This result is in a close
connection with the result of our
paper \cite{GyoMo} where we have determined all the linear bijections of
$C(X)$ (the algebra of all continuous complex valued functions on the
first countable compact Hausdorff space $X$) which preserve the diameter
of the range of functions. In fact, in $C(X)$ the spectrum of an
element $f$ is exactly its range. As the result in \cite{GyoMo} seems to
attract considerable interest among some researchers in the field of
function algebras, and there is so much interest in preserver problems
on operator algebras which concern the spectrum, we would like to pose
the following open problem.

\smallskip
{\it Problem.} Determine all the bijective linear transformations on
$B(H)$, the algebra of all bounded linear operators on the Hilbert space
$H$, which preserve the diameter of the spectrum.
\smallskip

Observe that our result Theorem~\ref{T:variance2} solves the
corresponding problem for $B_s(H)$. Regarding the mentioned facts,
we believe that this
is a prosperous and quite deep problem which deserves some attention.

\section{Acknowledgements}

The authors are very grateful to the referee for his/her
important observations and remarks that helped to improve the
presentation of the paper.

% Bibliography
\bibliographystyle{amsplain}

\end{document}